\pdfoutput=1
\documentclass[a4paper,conference]{IEEEtran}
\usepackage{graphicx,todonotes}
\usepackage{multirow}
\usepackage{amsmath,amsfonts,amssymb,amsthm}
\usepackage{cite}

\usepackage{todonotes}

\usepackage{caption} 
\usepackage{subcaption}
\usepackage[obeyspaces]{url}
\PassOptionsToPackage{obeyspaces}{url}
\usepackage{hyperref}

\begin{document}
	\bstctlcite{MyBSTcontrol}
	\title{Assessing Energy Storage Requirements\\Based on Accepted Risks}

\author{\IEEEauthorblockN{Michael P. Evans}
		\IEEEauthorblockA{Department of Electrical and Electronic Engineering\\
			Imperial College London, UK\\
			m.evans16@imperial.ac.uk}
		\and
		\IEEEauthorblockN{Simon H. Tindemans}
		\IEEEauthorblockA{Department of Electrical Sustainable Energy\\
			TU Delft, Netherlands \\
			s.h.tindemans@tudelft.nl}
	}
	
	\maketitle
    
\begin{abstract}
   This paper presents a framework for deriving the storage capacity that an electricity system requires in order to satisfy a chosen risk appetite. The framework takes as inputs user-defined event categories, parameterised by peak power-not-served, acceptable number of events per year and permitted probability of exceeding these constraints, and returns as an output the total capacity of storage that is needed. For increased model accuracy, our methodology incorporates multiple nodes with limited transfer capacities, and we provide a foresight-free dispatch policy for application to this setting. Finally, we demonstrate the chance-constrained capacity determination via application to a model of the British network.
\end{abstract}

\vspace{2mm}

\section{Introduction}
\label{sec:intro}
As electricity networks are decarbonised, an increasing proportion of the supply mix is provided by intermittent renewable sources. This leads to difficulty in balancing supply and demand. Energy storage offers a promising means of addressing this challenge, whereby fluctuations are smoothed by charging at times of oversupply and discharging at times of shortfall. This naturally leads to the question of how much storage a national electricity system requires. Due to economic constraints, it is not satisfactory to simply over-procure storage capacity by a large margin, but any level that has practical limits is necessarily accompanied by a risk of failure to balance the system. This is exemplified in the British (GB) network by the outage of 2019, during which there was simply not sufficient backup capacity available to meet the requirement of a rare event \cite{9Aug}.

In this work, we provide a framework for mapping risk appetites to storage capacity requirements. The user is able to define categories of shortfall event, and then specify how often they are willing to accept each type of event occurring. Chance-constrained specification of the form presented in \cite{Evans2019} is used to determine the minimum capacity of energy storage required. For the purposes of description, we focus on the GB system, but our approach remains general. Rather than use a static grid supply mix, we consider the evolution over time as predicted by National Grid, the GB TSO, in a scenario intended to meet the UK's obligations under the Paris Agreement \cite{FES}. 

Multiple authors have previously considered the question of the storage required by a system undergoing decarbonisation. Indeed, a summary of the results from multiple studies, including a fitted trend-line of storage capacity as a function of percentage renewable energy supply, was presented in \cite{ZERRAHN2018}, and then extended in \cite{Drax} to include National Grid outputs. The underlying methodology (see, for example, \cite{VICTORIA2019,CEBULLA2017,ZERRAHN2018, SCHOLZ2017,Macdonald2016}) involved minimising total system cost subject to specific decarbonisation conditions. However, this included the constraint that supply must meet demand in all locations and at all times. We extend this approach to include chance constraints, so that system planners willing to accept some risk of failure are able to perform a comparable optimisation task.

We design our event criteria based on prevalent generation adequacy metrics. Perhaps the two most common are expected energy-not-served (EENS) and loss-of-load expectation (LOLE). While the GB network has historically been designed based on LOLE, the inclusion of storage has led to the additional use of EENS for storage calculations, as a more representative encoding of event severity \cite{zachary2019}. It is for this reason that we focus our attention on the severity of shortfalls, in place of simply whether or not they occur. For the purposes of this work, we choose to consider peak power-not-served (PPNS) (i.e. maximum instantaneous shortfall across an event) under a specified dispatch strategy, but it would be straightforward to substitute this for an energy-related metric.

Our approach uses Monte Carlo simulation to sample combinations of supply and demand, each resulting in a shortfall trace. Integral to our analysis is then the simulated application of storage units to such a time-series. In \cite{Evans2017,Evans2018a}, a foresight-free dispatch policy was presented, in continuous and discrete time respectively, and shown to minimise the energy-not-served under pure discharging operation. An extension for recharging was also given in \cite{Evans2018a}. However, these policies apply to the case in which all storage units contribute towards meeting an aggregate shortfall at a single node; or equivalently a multi-node network with no limits on flows. Here, we are interested in the case in which transfer capacities are constrained. We present a quadratic programming equivalent to the round-trip policy of \cite{Evans2018a}, under the restriction that each device has unity efficiency and symmetric ratings for charging and discharging, that also allows for cross-charging among units. In addition to the current application, this could also be utilised in Model Predictive Control (MPC) dispatch frameworks.

\section{Problem description}
\label{sec:prob form}
\subsection{Mathematical description}
We denote by $n$ the number of nodes on the considered network, and utilise subscripts to denote node identifiers. At each node, our analysis is conducted based on annual supply and demand traces of sample period $\Delta t$ h, so that each trace has length $K\doteq 8760/\Delta t$. We consider a single nodal demand trace $D_i\colon\{1,...,K\}\mapsto[0,+\infty)$, but $m$ supply traces corresponding to different generator classes, each denoted $G_{i,j}\colon\{1,...,K\}\mapsto[0,+\infty)$. All of these traces will later be generated via stochastic processes. The nodal net shortfall $S_i[\cdot]$ can then be calculated according to
 \begin{equation}
 S_i[k]=D_i[k]-\sum_{j=1}^m G_{i,j}[k],\ k=1,...,K. 
 \label{eq:shortfall}
 \end{equation}

 In this work, we assume that the system planner is only able to apply storage to counteract shortfalls, and for the purposes of description restrict ourselves to a single storage unit per node. Note that this is without loss of generality, as co-located units can simply be connected by limitless lines. We denote the extractable energy stored at each node and sample instant as $e_i[k]$, which is subject to the assumed physical constraint $0\leq e_i[k]\leq \overline{e}_i$. We choose as our control input $p_i[k]$, the power extracted from the storage device, measured externally so as to take into account any inefficiencies present during discharging operation. This then leads to the following integrator dynamics:
\begin{equation}
e_i[k+1]=\left\{
\begin{array}{@{}ll@{}}
e_i[k]-p_i[k]\Delta t , & \text{if}\ p_i[k]\geq0 \\
e_i[k]-\eta p_i[k]\Delta t , & \text{otherwise},
\end{array}\right.
\label{eq:dynamics}
\end{equation}
in which $\eta_i\in[0,1]$ denotes the round-trip efficiency of the unit. We denote the power limits as $p_i[k]\in[-\underline{p}_i,\overline{p}_i]$.
Once storage has been dispatched, we refer to the adjusted shortfall trace as a \textit{resultant trace}, and each continuous period of remaining positive shortfall as a \textit{shortfall event}.

\subsection{Encoding of risk appetites}
\label{sec:risk appetites}
In this paper, we provide a means through which user-defined risk appetites can be mapped to storage requirements. We characterise the severity of an event based on PPNS, and rather than a binary check on whether a single threshold is exceeded, we allow the user to define $q$ event categories that are assessed independently. 

An example of $q=4$ category definitions can be found in Table~\ref{tab:event_classes}. For the $j$th event class, the user sets not only the PPNS threshold above which an event is registered, but also an allowed number of occurrences per year (total across all nodes), $a^{(j)}$, and an acceptable probability of failing to meet this tally, $c^{(j)}$. We define the counting operator $\Lambda^{(j)}(\cdot)$, which maps an annual resultant trace to an integer tally of the number of qualifying shortfall events. We then form the following chance-constrained optimisation problem, to determine the total storage power capacity required:
\begin{equation}
	    \begin{aligned}
	    &\min \overline{P}^{(j)} \\
	    \textnormal{s.t. } &\ \overline{p}_i=\rho_i\overline{P}^{(j)} \ \ \forall i \\
	    \textnormal{Pr}\ &\big[\sum_{i=1}^n\Lambda^{(j)}(S_i[\cdot]-p_i[ \cdot,\overline{p}_i])\leq a^{(j)}\big]\geq 1-c^{(j)},
	    \end{aligned}
	    \tag{$\star$}
	    \label{eq:optimisation problem}
\end{equation}
in which Pr$[\cdot]$ is the probability operator. The parameter $\rho_i$ denotes the (assumed fixed) proportion of total capacity allocated to the $i$th node, and the energy capacity is assumed to scale with the output power (i.e. the storage duration is assumed fixed), so that the battery capabilities are parameterised by the single parameter $\overline{P}^{(j)}$. In words, this optimisation problem aims to minimise the total storage capacity, subject to a chance constraint on the event tally across all nodes exceeding its allowed limit. Note that it is assumed that a given dispatch strategy is used (see next section). We repeat this full process for all $q$ event categories, and can additionally collapse the output into a single scalar value by computing $\max\{\overline{P}^{(j)},\ j=1,...,q\}$.

\section{Dispatch of storage under network constraints}
As mentioned in the \nameref{sec:intro}, we extend the dispatch policy of \cite{Evans2018a} to incorporate transfer capacity constraints. In order to allocate dispatch of the storage fleet at each time instant, a quadratic program is composed, in which the constraints correspond to the power and energy limits of each storage unit, as well as transfer ratings. This is then solved for each time step and an update to the energy performed according to the dynamics \eqref{eq:dynamics}. For clarity of notation we omit the explicit time dependence in the following. The objective function of the quadratic program consists of two sets of terms that describe a load shedding cost and a reduction in value of stored energy per node, each of which is described below.

\subsubsection{Avoidance of load shedding}
Our objective is to minimise the unserved demand after storage, which we denote as $s_i$ for the $i$th node and forms a decision variable. We assign a marginal cost that is linear in the proportion of demand that is not served, i.e.
\begin{equation}
    \alpha+\beta\frac{s_i}{D_i},
    \label{eq:shortfall objective}
\end{equation}
for positive $\alpha$ and $\beta$, leading to a total cost per time step of
\begin{equation}
    \frac{\beta\Delta t}{2D_i}s_i^2+(\alpha\Delta t) s_i.
\end{equation}

\subsubsection{Storage dispatch}
In the interest of simplicity, we assign the units full efficiency. As our modelling incorporates all losses into recharging operation, this should have a negligible effect when units can usually fully recharge between shortfall events; see \cite{Evans2018a} for further discussion on this point. We also set $\underline{p}=\overline{p}$.

The compound round-trip policy of \cite{Evans2018a} involves instantaneously targeting an equal value of \textit{time-to-go}, the remaining time for which each device can discharge at its maximum rate, across the fleet. During discharging operation, this involves ordering devices by descending time-to-go and allocating at maximum rate until the request is satisfied or all available devices are utilised \cite{Evans2017}. During recharging operation, this is inverted: devices are ordered by increasing time-to-go and allocated at maximum recharge rate \cite{Evans2018a}. By targeting an even time-to-go across the fleet, the policy aims to keep devices non-depleted, and therefore available for utilisation, under pure discharging request signals for as long as possible. From this intuition, one can deduce that devices with low time-to-go are inherently valuable; that is, there is a lower opportunity cost to  dispatching high time-to-go devices than there is to dispatching low time-to-go units.

When composing our objective function, we aim to capture this value, in the form of a value function, so that the returned optimum will push the time-to-go of each device in the correct direction. We achieve this by assigning a marginal value per unit of energy of 
\begin{equation}
    \gamma-\frac{\delta}{\overline{p}_i}e_i,
    \label{eq: dispatch objective}
\end{equation}
for positive $\gamma$ and $\delta$ satifying $0<\delta<\overline{p}_i\gamma/\overline{e}_i$. Integrating from $e_i=0$ results in a total value of
\begin{equation}
    V_k=\gamma e_i-\frac{\delta}{2\overline{p}_i}e_i^2
\end{equation}
at the beginning of a time slice. 
A comparison with the value function after an allocation of $p_i$ over one sample period leads to a \emph{reduction} in value (i.e. implied cost) of
\begin{equation}
    V_k-V_{k+1}=(\gamma-\frac{\delta}{\overline{p}_i} e_i) p_i\Delta t+\frac{\delta}{2\overline{p}_i}(p_i\Delta t)^2.
\end{equation}
We attempt to minimise this expression via the decision variable $p_i$. In order to capture that shortfall minimisation takes precedence over correct dispatch of storage units, we also set parameters such that \eqref{eq:shortfall objective} is larger than \eqref{eq: dispatch objective} for all units and areas.

\section{Solution methodology}
\subsection{Network model}
We partition the GB network into three geographical regions: Scotland, Northern England and Rest of GB. We then model the network as a complete three-node graph, where each node corresponds to one region. Transfer capacities between nodes are given by the bulk power transfer limits between regions. A flow network formulation was used in this paper, but this could be replaced by a power flow formulation if desired.

\subsection{Chance-constrained optimisation}
\begin{figure}[htb]
    \centering
    \includegraphics[width=\columnwidth]{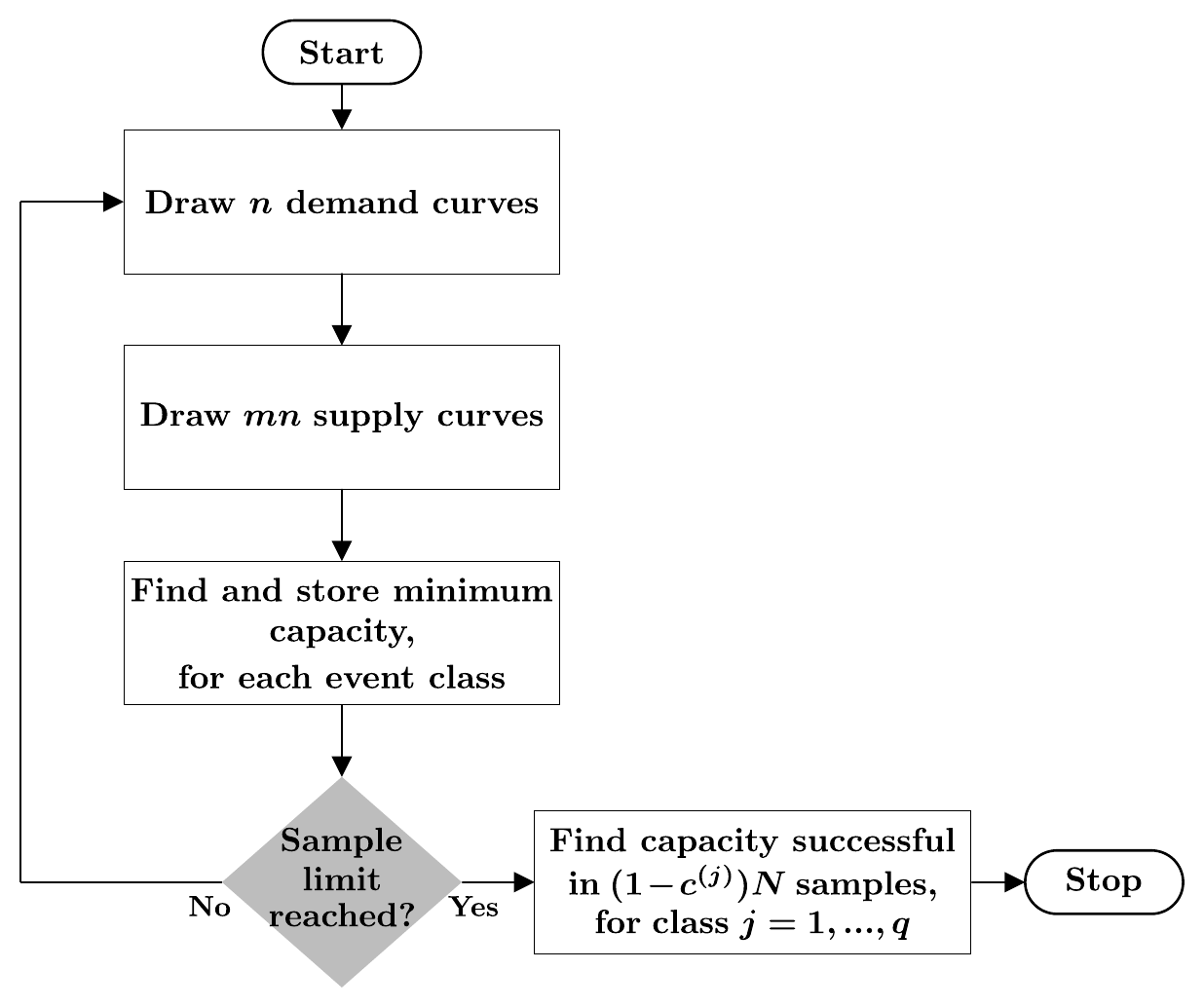}
    \caption{The Monte Carlo procedure}
    \label{fig:flow_chart}
\end{figure}

In this work, we utilise a comparable Monte Carlo procedure to that of \cite{Evans2019}. $N$ samples are composed by drawing demand and supply traces and computing, for each class of shortfall event, the boundary storage value at which the transition occurs from success to failure in adhering to the event criterion. A bisection routine is used to determine this value, with the output taken as the lower bound of a 1 GW interval. Once this routine has been repeated for all $N$ samples, we choose the value corresponding to success in $(1-c^{(j)})N$ samples, for event classes $j=1,...,q$, as can be seen in Figure~\ref{fig:flow_chart}.

\subsection{Trace generation}
The above procedure is repeated to return results for all future years up to 2050. In each case, annual supply and demand traces, of 1 h sample period, are generated according to the methodology of the following subsections. Once these have been generated, they are scaled to reflect National Grid's current prediction of the energy mix (under the `2 Degrees' scenario) \cite{FES} for the considered future year.

\subsubsection{Wind and solar}
We draw historical traces from the range 2000-2018, generated according to the reanalysis methodology of Staffell and Pfenninger \cite{Staffell2016,PFENNINGER2016}. One trace is composed per node, for each of offshore wind, onshore wind and solar. The location for each region is assigned to the largest (in terms of power rating) site in that region, which is used as the input coordinates to the data generator available at \cite{RN}. The output of this process is scaled to reflect the proportion of total generation of that class currently provided by the chosen region \cite{DUKES}. In order to account for geographical correlations, identical years are used for all solar traces, and identical years are used for all wind traces.

\subsubsection{Thermal plant, hydropower and marine}
The capacity for all thermal plant, including nuclear, is combined and then distributed among nodes according to current proportions \cite{DUKES}. The nodal total is then divided among generators of 500 MW rating. We repeat this generator assignment process for hydropower. We then compute traces for each generator via a 2-state Markov chain corresponding to up and down states, with a mean time between subsequent failure events of 2000 h. We set availability values based on the mean across relevant (Winter) data from Figure 3.1 of \cite{ECA}. As no marine generation is currently deployed on the system, we neglect any contribution for this source. However, this represents at most 2\% of predicted capacity moving forward in time.

\subsubsection{Interconnectors}
We assume the interconnectors that will be online from 2022 according to the market regulator \cite{ICs}, with the same proportion of total interconnection capacity allocated to each moving forward in time. We assign each to the nodes at which they are connected, then aggregate and divide the total capacity among 500 MW circuits that are able to fail independently. We model each of these circuits via the above 2-state Markov Chain procedure, this time with a 100 h mean time between subsequent failures (reflecting availability of power in the connected system). We assume that National Grid de-rating factors are a good representation of availability in this case, and so take the mean 2023/2024 value from \cite{CM}.

\subsubsection{Demand}
We draw aggregate demand traces from historical national demand traces, covering the period 2006-2018, which we scale to reflect predicted peak demand according to \cite{FES}. These data are provided by National Grid \cite{demand_data}, and converted from half-hourly to hourly averages. In order to disaggregate total demand, we divide this by region in proportion to total energy use for the same historical year \cite{BEIS_demand_by_region}.

\section{Results and discussion}
\subsection{Event categories}
 For our case study, we choose four event categories, as shown in Table~\ref{tab:event_classes} (recall that PPNS is the peak-power-not-served, and that an event is registered whenever the given value is exceeded). The boundaries between classes are defined based on 1, 3 and 5 multiples of the magnitude of the GB power cut on the 9 August 2019 \cite{9Aug}. It would be possible to extend our criteria to include considerations of recovery time, but we choose not to do so as it is difficult to predict the consequence of secondary effects. For example, following the 9 August event, poorly designed train post-fault protocols caused passengers to be stranded for significantly longer than it took to restore power \cite{9Aug}. 
\begin{table}[h]
    \renewcommand{\arraystretch}{1.5}
    \centering
    \begin{tabular}{|l|c|c|c|}
    \hline
     \textbf{Category}  & \textbf{PPNS (GW)} & \textbf{Number/year} & 
     \textbf{Exceedance probability} \\
     \hline \hline
     \textbf{Routine} & 0 & 3 & 1 year in 2 \\
     \hline
      \textbf{Mild} &  1.9 & 0 & 1 year in 5 \\
      \hline
      \textbf{Medium} & 5.7 & 0 & 1 year in 10 \\
      \hline
      \textbf{Severe} &
      9.5 & 0 & 1 year in 15 \\
      \hline
    \end{tabular}
    \caption{The categories of shortfall event considered.}
    \label{tab:event_classes}
\end{table}

\subsection{Three-node model}
\begin{figure}[htb]
    \centering
    \includegraphics[width=\columnwidth]{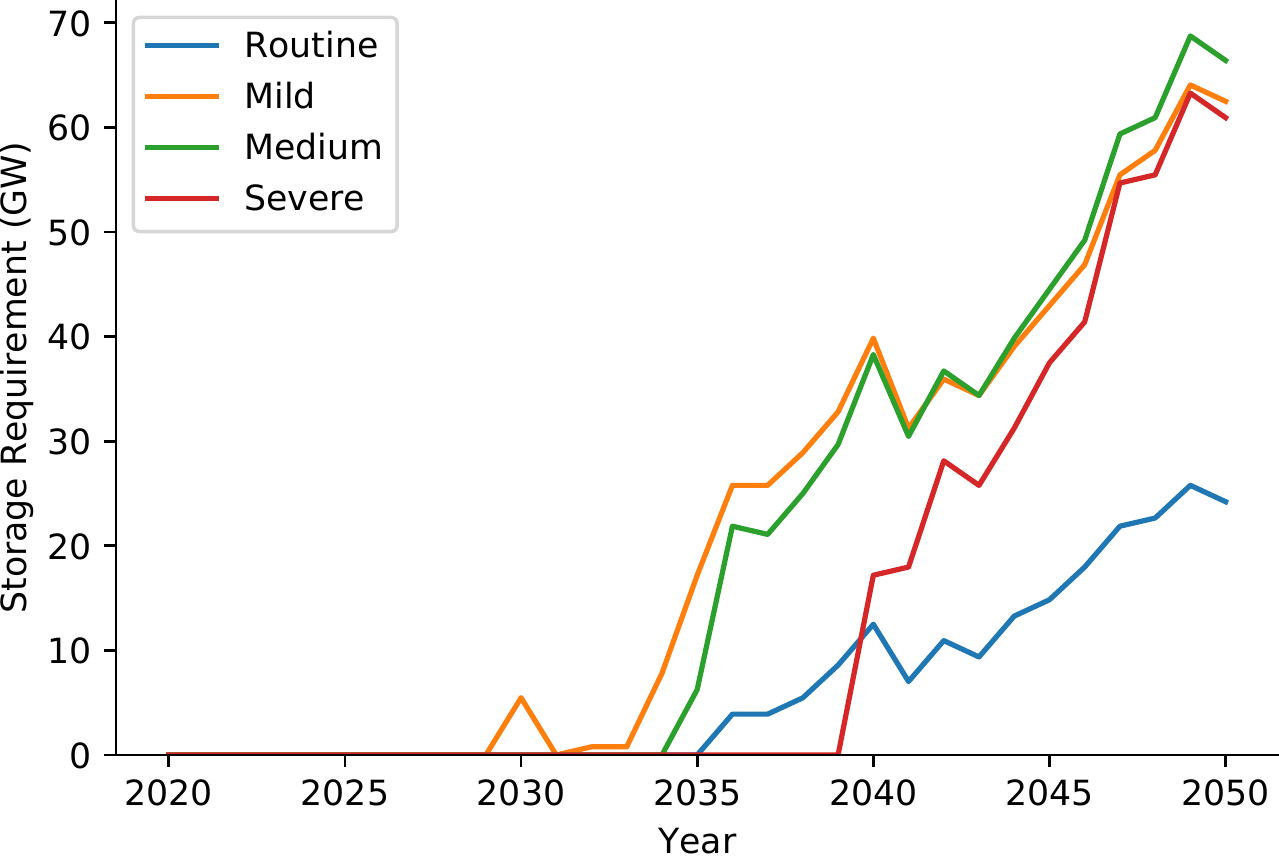}
    \caption{The progression of storage need over time, under the three-node model.}
    \label{fig:3_node_by_year}
\end{figure}

We perform the full study detailed above, across $N=10^4$ samples, with the total storage divided among nodes in proportion to demand. The observed progression of storage need over time can be seen in Figure~\ref{fig:3_node_by_year}. These results would need to be verified, for example by running the procedure for a longer period of time to ensure convergence to true values, and refinements made to modelling assumptions. The results do, however, demonstrate the potential impact of the framework that we present.

\subsection{Single-node model}
\begin{figure}[htb]
    \centering
    \includegraphics[width=\columnwidth]{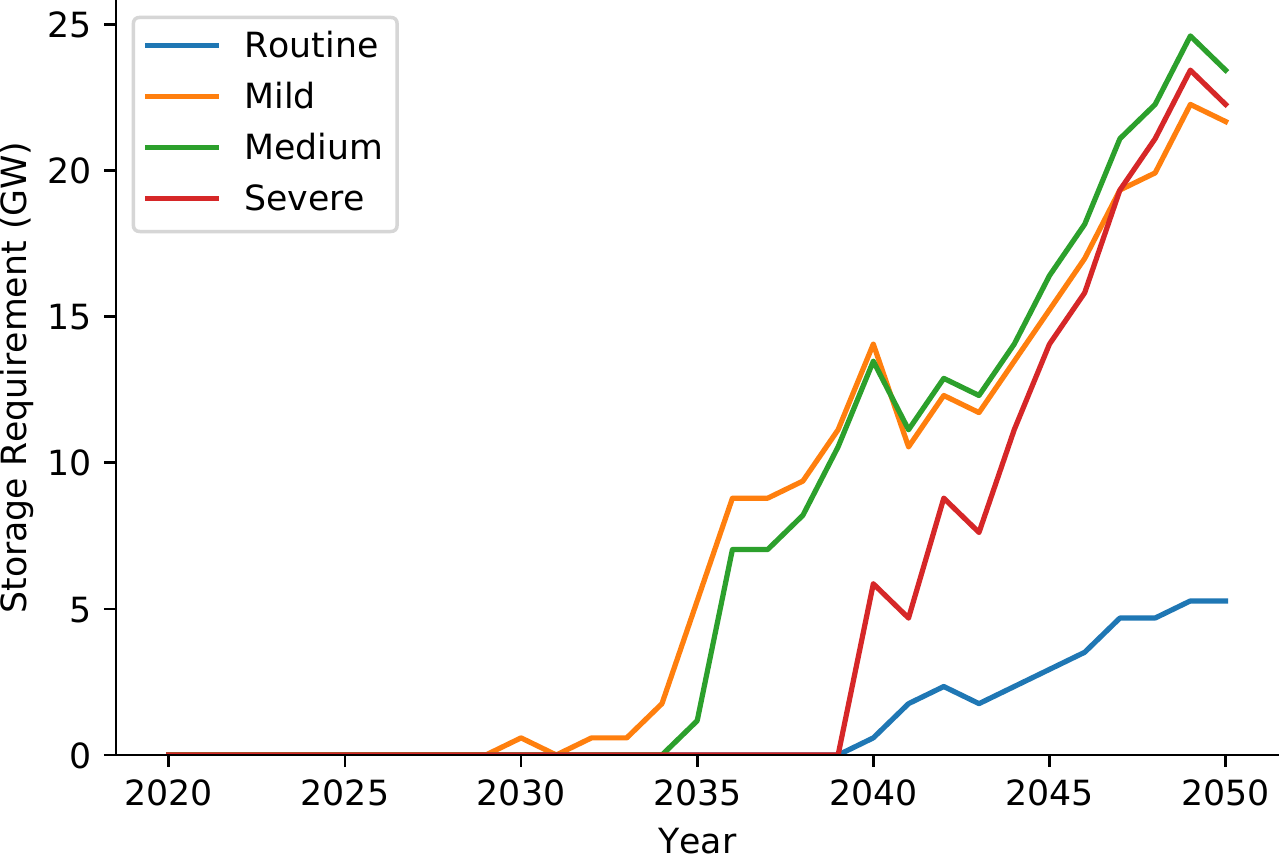}
    \caption{The progression of storage need over time, under the single-node model.}
    \label{fig:1_node_by_year}
\end{figure}

For comparison, we also perform the analysis for a single-node model. This is equivalent to the three-node model under the conditions that transfer capacities are not limited and with the applied storage composed of a single unit. The results of this investigation can be seen in Figure~\ref{fig:1_node_by_year}.

\subsection{Discussion}
Firstly, we compare the results from the three-node and one-node models as follows. We note that the trends are very similar in the two cases, including the ordering of event category requirements over time: in both cases, the Mild curve dominates up to approximately 2040, after which point it is the Medium curve that sets the maximum storage requirement. For all event categories and at all times, the 1-node results set lower bounds on the three-node results. This is as to be expected, since the combining of storage units and neglecting of transfer capacity limits should allow for better cross-compensation among nodes. The lack of optimisation in the distribution of total storage capacity among nodes will likely mean that the returned values for the three-node model are larger than the true storage needs for a three-node network. Conversely, limiting our analysis to three nodes and disregarding additional power flow constraints serves to reduce the storage requirement. Further investigation would be needed to quantify these effects for a particular case. 

The multi-node storage dispatch policy presented in this paper assumes that units prioritise addressing shortfalls over other objectives. This could reflect, for example, a system operator providing very large financial incentives for ancillary service participation. In cases where this compensation was in fact comparable to the potential earnings of the storage operator under alternative service routes, this modelling assumption might no longer hold. In order to roughly account for this, one could scale the resulting storage requirement by the inverse of the expected availability at the start of a shortfall event. For example, if the operator expected to have a 50\% state-of-charge at the start of each event, then the requirement would be doubled.

The exact storage requirement returned by our framework depends on the parameters of the model, for example: the choice of storage duration and conventional generation availability values; the distribution of storage among nodes; and the event category definitions. It would be up to the user of the framework, most likely a system planner, to set representative parameters and further investigate sensitivity to these choices.

For the purposes of this study, we have only considered total storage capacity as a decision variable, with all other capacities set to a constant value. A techno-economic optimisation that also considered changing the generation portfolio, or included flexible demand, might arrive at novel solutions to tackle the risks considered.

\section{Conclusion and Future work}
This paper has presented a framework for the chance-constrained determination of electricity system storage requirements, demonstrated through application to the GB network. The framework is capable of encoding multiple event categories, defined based on PPNS, allowed annual frequency and probability of exceeding this number of occurrences per year. This work has also introduced a foresight-free multi-node dispatch strategy for energy storage units that reduces to the policy of \cite{Evans2018a} for the single-node case. 

In future work, the authors intend to extend this model to include probabilistic transfer capacities and the effect of power flow constraints. They also plan to derive a theoretical basis for the value function used in dispatch optimisation, and in addition to extend this for cases where devices have imperfect efficiencies and asymmetric power ratings.

\section{Acknowledgements}
MP Evans' contribution to this project was funded by Eelpower Ltd.

\allowdisplaybreaks
\bibliographystyle{IEEEtran}
\bibliography{bib_all,bib_manual}
\end{document}